\documentclass[10pt]{article}
\usepackage{amssymb,amscd,amsmath,amsthm}
\usepackage{latexsym,amstext}
\usepackage{latexsym,amstext}
\usepackage{color}
\usepackage{graphicx}
\usepackage{epstopdf}

\begin{document}

\title{\bf An efficient implementation of the Hill-Harmonic Balance method to obtain  Floquet  exponents and solutions for homogeneous linear periodic differential equations}

\author{M. Gadella$^1$, L.P. Lara$^2$}

\maketitle

{$^1$ Departmento de F\'isica Te\'orica, At\'omica y Optica. University of Valladolid. Paseo Bel\'en 7, 47011, Valladolid, Spain, manuelgadella1@gmail.com.\\[2ex]
$^2$ Departamento de F\'isica, FCEIA, UNR, Av. Pellegini 250, 2000 Rosario, Argentina, lplara2014@gmail.com}

\begin{abstract}

We propose an implementation of a method based on Fourier analysis to obtain the Floquet characteristic exponents for periodic homogeneous linear systems, which  shows a high precision. This implementation uses  a variational principle to find the correct Floquet exponents among the solutions of an algebraic equation.   Once we have these Floquet exponents, we determine explicit approximated solutions. We test our results on systems for which exact solutions are known to verify the accuracy of our method.   Using the equivalent linear system, we also study approximate solutions for homogeneous linear equations with periodic coefficients.

\end{abstract}

\maketitle

\section{Introduction}

Linear periodic differential equations and systems of equations have an enormous presence in theoretical physics and engineering: harmonic oscillator, little oscillations, vibrations etc. However, there is a limited number of them that can be exactly solvable. In most cases, numerical methods are the only available. In order to obtain solutions of linear systems with periodic coefficients, one obtain the so called Floquet or characteristic exponents that determine a fundamental matrix for the system. This determination is usually given by a numerical approximation. In the case of linear equations with periodic non-constant coefficients, we always have the possibility of constructing the associated linear system, where the coefficients are again periodic, and then solving the system by means of the Floquet exponents.  

One of the most popular procedures to give approximate solutions for linear differential equations with periodic coefficients uses truncated Fourier series, whose coefficients are determined by the widely used  Harmonic Balance method \cite{BEL,CHEN,SAL,COC}. A modification of this method, which is particularly suitable for the Mathieu equation and other Hill type equations has been proposed in \cite{GGL}. 

Precisely the Harmonic Balance method has been used as an intermediate tool in order to obtain an approximation of the Floquet exponents for linear systems with periodic coefficients \cite{LAZ}. 

In the present paper, we introduce a modification of the procedure in \cite{LAZ} including a variational principle which gives the Floquet exponents as the critical values of this variational principle. This is easy to use and provides a great accuracy for the Floquet exponents and solutions as an added value. This efficiency as well as the accuracy of our modification has been tested in specific examples (like for instance in the Mathieu equation), as is well known the difficulty to establish error bounds for Fourier series based approximation methods. This method is primarily targeted to obtain the Floquet exponents of linear periodic systems,  although the application to obtain analytic algebraic approximate solutions of linear differential equations with periodic coefficients is then straightforward. 

This consideration of the Floquet coefficients as critical points of a variational problem is what makes our point of view different of previous methods including those in \cite{LAZ}. 

Before a description of our method and for the benefit of the reader, let us begin with an account of some important and well known results which are relevant in our presentation \cite{CL}. 

Let $A(t)$ be an $n\times n$ real matrix with continuous entries on the variable $t$ and ${\bf x}(t)\in\mathbb R^n$ for each value of $t$. In addition, all these entries are periodic with the same period $T$, so that $A(t+T)=A(t)$ for all $t$.  Let us consider a linear system of the form

\begin{equation}\label{1}
\dot{\bf x}(t)=A(t){\bf x}(t)\,,
\end{equation}
where the dot means derivative with respect to $t$. The Floquet theory which refers to this type of systems is well known \cite{CL,CI,VER}. In the sequel, we recall some interesting well known facts which are useful in our discussion \cite{CL}:

\medskip

i.) If $\Phi(t)$ is a fundamental matrix of (\ref{1}), $\Phi(t+T)$ is again a fundamental matrix. As is the case for any pair of two fundamental matrices, there must be a constant invertible matrix $C$ such that 

\begin{equation}\label{2}
\Phi(t+T)=\Phi(t)\,C\,.
\end{equation}

Since $C$ is invertible, it must exist a $n\times n$ matrix $B$ such that

\begin{equation}\label{3}
C=e^{BT}\,,
\end{equation}
where $T$ is again the period for $A(t)$. 

\medskip

ii.)  Consider the matrix $P(t):=\Phi(t)\,e^{-Bt}$. Then, $P(t)$ is periodic with period $T$ and $P(t)$ is invertible.

\medskip

iii.) Let us consider the following new indeterminate ${\bf y}(t)$ as:

\begin{equation}\label{4}
{\bf y}(t)=P^{-1}(t)\,{\bf x}(t)\,.
\end{equation}

Then  and taking into account that $P(t)$ is invertible, we have that

\begin{equation}\label{5}
\dot{\bf y}=B\,{\bf y}\,.
\end{equation}

Thus, system (\ref{1}) is equivalent to a system with constant coefficients. We shall recall in a moment on the importance of the matrix $B$. 

\medskip

iv.) Then, if we define an initial condition ${\bf x}(t_0):=P(t_0){\bf y}(t_0)$ and taking into account that that the solution of (\ref{7}) satisfying the initial condition ${\bf y}(t_0)$ is given by

\begin{equation}\label{6}
{\bf y}(t)=e^{(t-t_0)B}\, {\bf y}(t_0)\,,
\end{equation}
we have that

\begin{equation}\label{7}
{\bf x}(t)=P(t)\,e^{(t-t_0)B}\, P^{-1}(t_0)\,{\bf x}(t_0)\,.
\end{equation}

Let us choose as $\mathbf y(t_0)$ and eigenvector $\mathbf y_0$  of $B$ with eigenvalue $\lambda_0$, i.e., $B\,\mathbf y_0=\lambda_0\,\mathbf y_0$, where $\lambda_0$ is any of the eigevalues of $B$.  Then if $\mathbf x_0:=P(t_0)\,\mathbf y(t_0)=P(t_0)\,\mathbf y_0$, one has

\begin{eqnarray}\label{8}
\mathbf y(t)=e^{(t-t_0)\lambda_0}\,\mathbf y_0  \Longrightarrow \mathbf x(t)=P(t)\, e^{(t-t_0)\lambda_0}\,P^{-1}(t_0)\,\mathbf x_0\,,
\end{eqnarray}
expression which may be written as

\begin{equation}\label{9}
\mathbf x(t)=\eta(t)\,e^{(t-t_0)\lambda_0}\,,\qquad {\rm with}  \qquad \eta(t)=P(t)\,P^{-1}(t_0)\,\mathbf x_0\,.
\end{equation}

Since $P(t)$ is periodic with period $T$, equation (\ref{9}) shows that $\eta(t)$ is also periodic with period $T$. 

In summary, we can obtain particular solutions of (\ref{1}) if we can determine the eigenvalues of the matrix $B$. This eigenvalues are usually called {\it Floquet characteristic exponents} or {\it Floquet exponents} or simply {\it characteristic exponents}. We shall keep this terminology along our manuscript.  There are not general analytic methods to obtain these characteristic exponents and, hence, numerical methods for their determination are in order. 

In the present article, we propose an analytic approximate method in order to obtain the characteristic coefficients, with the following order:  In Section 2, we give a standard method to obtain the Floquet characteristic coefficients, important for a comparison with our proposed method. We introduce our analytic algebraic approximation method in Section 3, more precisely on 3.1, where we propose the variational principle to obtain the Floquet critical exponents. Section 4 is devoted to a test  using the Mathieu equation.  In Section 5, we test the method when $A(t)$ as in (\ref{1}) is $2\times 2$ and  commutes with its integral.  In this case, we give an explicit expression for the fundamental matrix $\Phi(t)$ of the system. In addition, we provide of an interpretation of the eigenvalues of the average of the matrix $A(t)$ over a period: the eigenvalues of the resulting matrix are the characteristic exponents. We test our results with a couple of examples, one being the Marcus-Yamabe equation. We close the paper with some concluding remarks. 

\section{Determination of the characteristic exponents: standard method}

Let us go back to Equation (\ref{1}), in which  $A(t)$ is periodic with period $T$.  As initial values, we may choose any of the vectors of the canonical basis in $\mathbb R^n$, i.e., those vectors with all components equal to zero except the i-th component which is equal to one.  Once we have chosen an initial value, a numerical  integration such as a fourth order Runge-Kutta  \cite{KIN} permits us to obtain $n$ discrete linearly independent solutions on the finite interval $(0,T)$, where $T$ is the period. Then, by using interpolation, for instance with splines, we obtain an approximate continuous solution.  Using the initial conditions, we obtain $n$ approximate linearly independent solutions $X_1(t),X_2(t),\dots,X_n(t)$, whose columns determine an approximate fundamental matrix $\Phi(t)$. This procedure is rather simple for $n=2$, which will be our case. 

After (\ref{2}) and (\ref{3}), we readily obtain

\begin{equation}\label{10}
C:=\exp\{BT\}=\Phi^{-1}(0)\,\Phi(T)\,.
\end{equation}

The relation between the eigenvalues $\Delta_i$ of $C$ and the characteristic coefficients $\lambda_i$ is well known:

\begin{equation}\label{11}
\lambda_i=\frac 1T\,\log \Delta_i\,,\qquad i=1,2,\dots,n\,.
\end{equation}

Thus, we have determined the characteristic coefficients and the numerical solution $X(t)$. We have to take into account that the imaginary part of the characteristic coefficients is not uniquely determined since:

\begin{equation}\label{12}
\Delta_i=\exp\{ \lambda_i+2\pi i/T\}T=e^{\lambda_i T}\,.
\end{equation}

Our choice will always fix this imaginary part in such a way that the exponent coincides with $\lambda_i T$, being $\lambda_i$ an eigenvalue of $B$.  

The objective of the present article is to show that a good approximation on the characteristic coefficients may be  obtained through an algebraic analytic approximation based on Fourier analysis. 

\section{Approximated analytic solution}

The relation between first order linear systems of the form (\ref{1}) and linear equations of order $n$ is well known \cite{CL}. With this idea in mind, let us illustrate our method with  second order linear differential equations of the form:

\begin{equation}\label{13}
\ddot y(t)+a(t)\dot y(t)+b(t) y(t)=0\,,
\end{equation}
where $a(t)$ and $b(t)$ are periodic functions with respective periods $T_a$ and $T_b$, which are not arbitrary, since we have to impose the condition that the ratio $T_a/T_b$ be rational. In addition, $a(t)$ is continuously differentiable and $b(t)$ continuous. The linear system equivalent to (\ref{15}) is ($z_1(t):=y(t)$, $z_2(t)=\dot y(t)$)

\begin{equation}\label{14}
\left( \begin{array}{c}  \dot z_1(t)\\[2ex] \dot z_2(t)   \end{array}  \right)  = \left( \begin{array}{cc}  0 & 1\\[2ex] -b(t) & -a(t)  \end{array}   \right)  \left( \begin{array}{c}  z_1(t)\\[2ex]  z_2(t)   \end{array}  \right) \Longleftrightarrow \dot {\mathbf z}(t)=A(t)\,\mathbf z(t)\,.
\end{equation}

With the change

\begin{equation}\label{15}
y(t)=x(t)\,\exp\left\{ -\frac12\,\int a(t)\,dt    \right\}\,,
\end{equation}
equation (\ref{15}) yields to

\begin{equation}\label{16}
\ddot x(t)+f(t)\,x(t)=0\,,
\end{equation}
with

\begin{equation}\label{17}
f(t)=b(t)-\frac12\,a'(t)-\frac14\,a^2(t)\,.
\end{equation}

The function $f(t)$ is continuous and periodic with a period $T=\max(T_a,T_b)$.  By the Floquet characteristic exponents, or just characteristic exponents, of  (\ref{13}), we mean the characteristic exponents of the associated system (\ref{14}). Analogously, the characteristic exponents of (\ref{16}) are the characteristic exponents of its related system. 

Let us list in the sequel some of the properties of equation (\ref{13}):

\begin{itemize}

\item{Assume that

\begin{equation}\label{18}
f(t)>0 \qquad {\rm and} \qquad T\int_0^T f(t)\,dt\le 4\,.
\end{equation}

Then, it has been proven in \cite{CI} that all solutions are bounded. Consequently, the characteristic exponents of (\ref{13}) do not have positive real part.  }

\item{If $f(t)<0$, let us multiply (\ref{13}) by $y(t)$ and integrate by parts. Then, we have 

\begin{equation}\label{19}
\frac{d}{dt}\,y^2(t)=\int [y'(t)]^2\,dt-\int f(t)\,y^2(t)\,dt>0\,.
\end{equation}

Since $y^2(t)\le 0$, we note that for large values on the variable $t$, $t\longmapsto\infty$, the solution $y(t)$ is not bounded.

Consequently, the characteristic coefficients must have a positive real part. 
}

\item{Let us go to equation (\ref{16}). It can be proven \cite{CI,VER} that  the sum of its characteristic exponents is equal to zero. }

\end{itemize}

\subsection{The method}

Consider equation (\ref{16}) and assume that $\lambda$ is one of its characteristic exponent. Choosing for simplicity $t_0=0$, we go back to (\ref{16}) were  it was stated that for each Floquet exponent there is a solution of the type $x(t)=\eta(t)\,e^{t\lambda}$, where $\eta(t)$ is periodic with period $T$\footnote{We may assume that a basis of solutions is of this form, provided that $B$ be diagonalizable.}. The point is that $\eta(t)$ and $\lambda$ are unknown and our objective is to find an approximate expression for them.   Using this result in (\ref{16}), we obtain the following differential equation:

\begin{equation}\label{20}
\ddot \eta(t)+2\lambda\,\dot\eta(t)+(\lambda^2+f(t))\,\eta(t)=0\,.
\end{equation}

We have obtained a second order equation with a periodic coefficient $f(t)$ with period $T$. Let us span $\eta(t)$ into Fourier series and then truncate this series. The truncated solution $\eta_n(\lambda,t)$ is

\begin{equation}\label{21}
\eta_n(\lambda,t)= \frac{a_0}2+ \sum_{k=1}^n \{a_k\cos (k\omega t)+b_k\sin(k\omega t)\}\,,
\end{equation}
with $\omega:= 2\pi/T$. Now, $x_n(t)=z_n(t,\lambda)\exp(\lambda t)$. 

In order to determine the characteristic exponent $\lambda$, we propose the following strategy:

\medskip

i) First of all, we determine the coefficients $a_k$ and $b_k$ by means of the Harmonic Balance (HB) method \cite{SAL,CHEN,BEL}. In summary, we replace (\ref{21}) into (\ref{20}) so as to obtain a new Fourier polynomial. Since equation (\ref{20}) must be satisfied, coefficients for the harmonics in this Fourier polynomial must vanish. This  yields to an homogeneous linear algebraic system of  dimension $2n+1$, with indeterminates $a_0$,  $a_k$ and $b_k$ and $k=1,2,\dots,n$. In order to obtain non-trivial solutions, the determinant $\Delta$ of the matrix of the system of the coefficients must vanish.  Since (\ref{20}) is linear, this determinant is a polynomial on $\lambda$, so that

\begin{equation}\label{22}
\Delta(\lambda)=0
\end{equation}  
gives $\lambda$ in terms of $\omega$ and any other parameter appearing in (\ref{20}).  Although (\ref{22}) has at most $n$ roots, only two of them could be the characteristic exponents we are looking for. Moreover, it is not difficult to check that the coefficients $a_k$ and $b_k$,  $k=1,2,\dots,n$ are rational polynomial functions on $\lambda$. 

\medskip

ii) After we have completed the previous step, we shall determine the approximate value of $\lambda$ by a variational principle. Since the exact solution $x(t)$ satisfies 

\begin{equation}\label{23}
\int_0^T (\ddot x(t)+f(t)\,x(t))^*\,(\ddot x(t)+f(t)\,x(t))\,dt=0\,,
\end{equation}
where the star denotes complex conjugation, we propose that the {\it approximate characteristic exponent}, $\lambda_k$, we are searching for is a {\it critical point} (usually a minimum) of $E(\lambda)$ defined as:

\begin{equation}\label{24}
E(\lambda):= \int_0^T (\ddot x_n(t)+f(t)\,x_n(t))^*\,(\ddot x_n(t)+f(t)\,x_n(t))\,dt\,.
\end{equation}

Note that $\lambda$ may have an imaginary part and thus $x(t)$. This is the reason why we have to include a complex conjugation in (\ref{23}-\ref{24}). 

Once we have the Floquet characteristic exponents for the given equation, we determine the coefficients $a_0$, $a_k$ and $b_k$, $k=1,2,\dots,n$ for the truncated Fourier series that approximates the solution. 

Our variational principle is just an Ansatz, which should be confirmed by numerical experiments. This is the objective of the rest of the article. 

\section{Application I: The Mathieu equation.}

The Mathieu equation is a simple non-trivial equation with periodic coefficients which is well suitable as a laboratory in order to test the above ideas as shown by previous work of our group \cite{GGL}. The Mathieu equation has been largely studied,  as for instance in \cite{DNH,FP,MW,BC,ZRKBS,AS}. Let us write the Mathieu equation as

\begin{equation}\label{25}
\ddot x(t)+\omega^2(1-\alpha\,\cos t)\,x(t)=0\,.
\end{equation}

As is well known, two linearly independent solutions are

\begin{equation}\label{26}
x_1(t)=C\left( 4\omega^2,2\alpha\omega^2,\frac t2  \right)\,, \qquad x_2(t)=S\left( 4\omega^2,2\alpha\omega^2,\frac t2  \right)\,,
\end{equation}
where $C$ and $S$ stand for the Mathieu sine and cosine \cite{AS}. These are exact solutions, so that we can determine {\it exact} characteristic exponents just by constructing a fundamental matrix with them and, then, making use of equations (\ref{11}) and (\ref{12}), which in this case give the exact results. 

Now the objective is clear and is to compare the results obtained with our proposed variational method with the exact results that can be obtained as described above. In addition, we shall also compare both with those obtained following the lines introduced in Section 2. 

Before proceeding, a couple of comments. First of all, using (\ref{18}) we see that for $\omega<1/4$ and for all values of $\alpha$ the solutions are bounded. Also note that whenever $\lambda=ik$, being $k$ integer, the solution is periodic with period equal to $2\pi$. Finally, let us recall that the sum of the critical exponents is equal to zero, an interesting property to take into account as we test our results. 

Let us go back to the determinant (\ref{22}), that we write now as $\Delta(\lambda)\equiv \Delta_{\alpha,\omega}(\lambda)$, due to its dependence on all these three variables. In our case, it is an even polynomial of degree $2(2n+1)$. Furthermore, in all cases studied it is also an even polynomial on the variables $\alpha$ and $\omega$. As an example, let us take $n=2$, so that the polynomial on $\lambda$ has degree ten:

\begin{equation}\label{27}
\Delta_{\alpha,\omega}(\lambda)=\sum_{k=0,{\rm even}}^{10} c_k\,\lambda^k\,.
\end{equation}

In (\ref{27}) all odd coefficients vanish, while the even coefficients are given by:

\begin{eqnarray}
c_{0}    =  16\,\omega^{2}+\,8\,(-5+\alpha^{2})\,\omega^{4}+(33-14\,\alpha
^{2})\,\omega^{6} \alpha^{2} -\frac{1}{2}(20-14+\,\alpha^{4})\,\omega
^{8}-\nonumber\\[2ex]-8\,(1-\,\alpha^{2}+\frac{3}{16}\alpha^{4})\,\omega^{10}\,,\nonumber
\end{eqnarray}
\begin{eqnarray}
c_{2}    =  \frac{1}{16}(256+(560-96\alpha^{2})\omega^{4}+160(\alpha
^{2}-2)\omega^{6}+(80-48\alpha^{2}+3\alpha^{4})\omega^{8})\,,\nonumber
\end{eqnarray}
\begin{eqnarray}
c_{4}   =  40+35\omega^{2}+3\alpha^{2}\omega^{4}+(10-3\alpha^{2})\,  
\omega^{6}\,, \qquad \qquad \qquad \qquad\qquad\qquad\nonumber
\end{eqnarray}
\begin{eqnarray}
 c_{6}    =33+20\omega^{2}-(-10+\alpha^{2})\omega^{4}\,,\qquad\qquad\qquad\qquad\qquad\qquad\qquad\qquad\nonumber
\end{eqnarray}
\begin{eqnarray}
c_{8}   =5(2+\omega^{2})\,, \qquad\qquad\qquad\qquad\qquad\qquad\qquad\qquad\qquad\qquad\qquad\nonumber
\end{eqnarray}
\begin{equation}\label{28}
c_{10}   =304+80\,\omega^{2}\,. \qquad\qquad\qquad\qquad\qquad\qquad\qquad\qquad\qquad\qquad\qquad
\end{equation}

Then, using the Harmonic Balance method that, in this case, is a simple algebraic problem in which the equations that determine the coefficients $a_k$ and $b_k$ are homogeneous and starting with the initial condition $a_1=1$, we obtain for the first coefficients the following values:

\begin{eqnarray}\label{29}
b_{1}  &  = & (-\lambda\alpha^{2}\omega^{4}+2\lambda(-16\lambda^{2} -(-4+\lambda^{2}+\omega^{2})^{2}))/d_{1} \,, \nonumber\\[2ex]
b_{2}  &  = & -12\lambda\alpha\omega^{2}(-4+\lambda^{2}+\omega^{2})/d_{2}\,,
\end{eqnarray}
\begin{eqnarray}\label{30}
a_0 =\frac{1}{\alpha\omega^{2}}\{3(-2+\lambda^{2}+\omega
^{2})-\frac{1}{8\lambda}[(16+4\lambda^{4}-20\omega^{2} 
-(-4+\alpha^{2})\omega^{4}+\nonumber\\[2ex]+\lambda^{2}(-52+8\omega^{2}))(-\lambda\alpha^{2}\omega^{4}+2\lambda(-16\lambda^2-(-4+\lambda^{2}+\omega^{2})^{2})]/d_{1}\}\,, \nonumber\\[2ex]
a_1=1\,,\nonumber\\[2ex]
a_2= -\frac{1}{2}[\alpha\omega^{2}(-16-4\lambda^{4}+20\omega
^{2}+(-4+\alpha^{2})\omega^{4}+\lambda^{2}(52-8\omega^{2}))]/
d_2\,,
\end{eqnarray}
where,

\begin{eqnarray}\label{31}
d_{1}   =\frac{1}{4}\alpha^{2}\omega^{4}(-4+\lambda^{2}+\omega
^{2})+(-1+\lambda^{2}+\omega^{2})(-16\lambda^{2}-(-4+\lambda^{2}+\omega
^{2})^2\,, \nonumber\\[2ex]
d_2= 4\lambda^{6}+4\lambda^{4}(7+3\omega^{2})+\lambda^{2}(32-8\omega
^{2}-(-12+\alpha^{2})\omega^{4})-\qquad \nonumber\\[2ex]
-(-4+\omega^{2})(-16+20\omega^{2}+(-4+\lambda^{2})\omega^{4}))\,.\qquad
\end{eqnarray}

We may had obtained similar expressions for higher values of $n$, although they are increasingly complicated and do not provide of any new information. Once we have obtained the roots $\lambda_k$ of (\ref{22}), only two of them can be chosen to be the critical exponents. They are precisely those which minimize (\ref{24}). Once we have obtained the critical exponents, we readily determine an explicit approximated solution of (\ref{16}). 

As an example, let us choose $\alpha=0.5$, $\omega=1$ and $n=2$. We obtain the following approximate solution:

\begin{equation}\label{32}
x_A(t)= \exp\left(-\frac{1}{43}t\right)\left(\frac{267}{1069}+\cos t-\frac{4}{45}\cos2t-\frac{185}{84}\sin t+\frac{15}{83}\sin2t\right)\,.
\end{equation}

Note that the general exact solution has the form $x_e(t)=c_1\,x_1(t)+c_2\,x_2(t)$, where $x_i(t)$, $i=1,2$ are given in (\ref{26}). The constants $c_i$, $i=1,2$ should be determined through the initial conditions $x_e(0):=x_A(0)$ and $\dot x_e(0):=\dot x_A(0)$, where the dot represents derivative with respect to $t$. These values are obtained with the expression for $x_A(t)$ in (\ref{32}).  

We know the exact value of the characteristic exponents, which validates our comparison,  take one and denote it as $\lambda_e$. These characteristic exponents are $\lambda=\pm 1/43$. In order to compare the exact solution with the approximation given in (\ref{32}), it is natural to choose the exponent with minus sign, so that $\lambda_e=-1/43$. Thus, the exact solution has the form $x_e(t)=e^{\lambda_e t}\,\eta(t)$, see comments before (\ref{20}). Since we have determined already $x_e(t)$ through the above initial conditions, we know $\eta(t)$. Then, span $\eta(t)$ into Fourier series. We obtain an explicit expression of the form:

\begin{eqnarray}\label{33}
x_e(t)= \exp\left(-\frac{1}{43}t\right)\,\left(\frac{66}{265}+\frac{486}{487}\cos
t-\frac{4}{45}\cos2t+\right.\nonumber\\[2ex] \left.+\frac{1}{340}\cos3t-\frac{131}{60}\sin t+
\frac{7}{39}\sin2t-\frac{1}{176}\sin3t+...\right)\,. 
\end{eqnarray}

The coefficients in both solutions have been adjusted to a rational number with an error upper bound of $0.07\%$. The relative difference between approximate (\ref{32}) and exact (\ref{33}) solutions is at most less than $0.9\%$. This is certainly satisfactory. As expected, a higher value of $n$ gives a higher precision. For instance, take $n=3$, $\alpha=0.5$ and $\omega=1$. We have for the approximate and exact solution, respectively, the following results:

\begin{eqnarray}\label{34}
x_A(t)= \exp\left(-\frac{1}{43}t\right)\,\left(\frac{267}{1069}+\cos t-\frac{5}{56}\cos2t+\right. \nonumber\\[2ex] \left. +\frac{1}{346}\cos3t-\frac{376}{171}\sin t+  \frac{19}{105} \sin2t-\frac{1}{178}\sin3t\right)
\end{eqnarray}
and

\begin{eqnarray}\label{35}
x_e(t)= \exp\left(-\frac{1}{43}t\right)\,\left(\frac{285}{1141}+\cos t-\frac
{5}{56}\cos2t+\frac{1}{346}\cos3t- \right.\nonumber\\[2ex]  \left. -\frac{596}{271}\sin t+\frac{19}{105}  \sin2t-\frac{1}{178}\sin3t\right)+...
\end{eqnarray}

For $n=3$, $\alpha=1$ and $\omega=1$, we obtain analogously:

\begin{eqnarray}\label{36}
x_A(t)= \exp\left(-\frac{1}{10}t\right)\,\left(\frac{1}{2}+\cos t-\frac{1}{5} \cos2t+\frac{1}{69}\cos3t-\right. \nonumber\\[2ex] \left.- \frac{19}{9}\sin t+\frac{1}{3}\sin2t-\frac{1}{50}\sin3t\right)
\end{eqnarray}
and

\begin{eqnarray}\label{37}
x_e(t)= \exp\left(-\frac{1}{10}t\right)\,\left(\frac{1}{2}+\cos t-\frac{1}{5}\cos2t+\frac{1}{69}\cos3t-\right. \nonumber\\[2ex]\left.-\frac{19}{9}\sin t+\frac{1}{3}\sin2t-\frac{1}{51}\sin3t\right)+...
\end{eqnarray}

It is important to stress that for $n=2$ and $n=3$, we have used different initial conditions so that the exact solutions (\ref{35}) and (\ref{36}) do not coincide. These initial conditions are given by the values of $x_A(t)$ and its first derivative at the origin. In particular, for $n=2$, we have $x(0)=1.16088$ and $\dot x(0)=-1.86793$. For $n=3$, we have $x(0)=1.16337$ and $\dot x(0)=-1.88083$. Since we have changed the initial conditions, we have changed the solution and therefore the critical exponents could be different, which is the case here.  

Observe that we have achieved a better precision. The conclusion is that the higher the harmonic number $n$ is the better accuracy is obtained. This result is quite satisfactory. 

\vskip1cm

$
\begin{array}
[c]{cccccc}
\alpha & 1/10 & 3/10 & 5/10 & 7/10 & 1\\[2ex]
\lambda_e & 9.31603\,10^{-4} & 8.37695\,10^{-3} & 2.32152\,10^{-2} &
4.52826\,10^{-2} & 9.10175\,10^{-2}\\[2ex]
\lambda_A & 9.31603\,10^{-4} & 8.37695\,10^{-3} & 2.32152\,10^{-2} &
4.52825\,10^{-2} & 9.10172\,10^{-2}\\[2ex]
\lambda_{num} & 9.31620\,10^{-4} & 8.37697\,10^{-3} & 2.32151\,10^{-2} &
4.52826\,10^{-2} & 9.10175\,10^{-2}\\[2ex]
S^{\,2} & 0. & 4.\,10^{-9} & 9.\,10^{-8} & 6.\,10^{-7} & 4.\,10^{-6}\\[2ex]
E(\lambda) & 2.10^{-10} & 4.10^{-7} & 6.10^{-6} & 7.10^{-5} & 5.10^{-4}
\end{array}
$

\vskip1cm

In Table 1, we compare the values of the approximate characteristic exponent given by our method, $\lambda_A$, the exact, $\lambda_e$, and the one detemined by the method sketched in Section 2, $\lambda_{num}$, for $n=3$, $\omega=1$ and different values of $\alpha$.  The precision of $x_A(t)$ is evaluated through the second moment

\begin{equation}\label{38}
S^2:=\frac 1T \int_0^T (x_e(t)-x_A(t))^*\,(x_e(t)-x_A(t))\,dt\,.
\end{equation}

\begin{figure}
\centering
\includegraphics[width=0.5\textwidth]{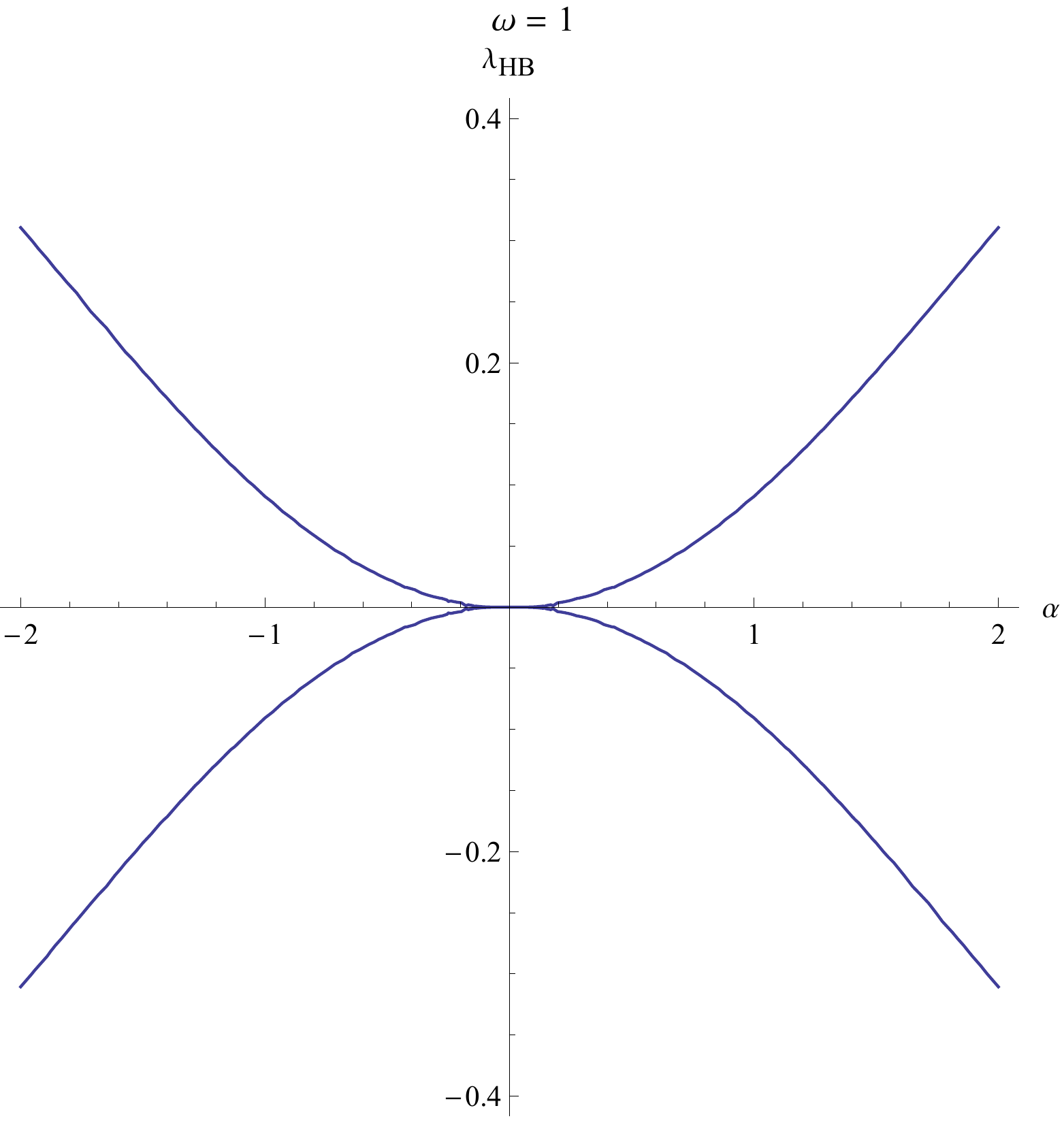}
\caption{Variation of $\lambda_A$ in terms of $\alpha$ for $n=2$ and $\omega=1$. 
\label{FIG1}}
\end{figure}

Finally, we include $E(\lambda)$ given in (\ref{24}), which measures the deviation of the solution $x_A(t)$ from the exact solution of the differential equation (\ref{16}). Since the sum of the exponents vanish, we refer only to one of them. Nevertheless, we must say that errors in the method which are always present in this kind of estimations, make the sum of both critical exponents not exactly equal to zero. An estimation with six digits of $\lambda_A$ exactly matches the exact result, while $\lambda_{num}$ has a minor discrepancy of three units in the last digit. Computational results have been preformed with the use of Mathematica, the CPU time being negligible.

In Figure 1, we plot the dependence of $\lambda_A$ with $\alpha$ for $n+2$ and $\omega=1$. Note that both solutions appear symmetric with respect to the abscise axis. Recall that there are always two solutions whose sum is equal to zero.

\section{Application II: Matrices commuting with their integral}

In \cite{TW}, the authors discuss the case in which $A(t)$ and the matrix obtained by integrating its entries with respect to the variable $t$ commute. We are going to explore how to apply our method when this is the case. It is noteworthy that we now have an explicit expression for the fundamental matrix $\Phi(t)$.   Let us go back to (\ref{1}) and  pose a result that has been proven in  \cite{TW}. Under the following hypothesis:

\smallskip

i.)  All entries of $A(t)$ in (\ref{1}) are integrable on the interval $[0,t]$. 

\smallskip

ii.) The matrix $A(t)$ fulfils the following commutation relation:  

\begin{equation}\label{39}
\left[  \int_0^t A(v)\,dv  \right] A(t)=A(t)  \left[  \int_0^t A(v)\,dv  \right]\,,
\end{equation}
where the integral of a matrix is the matrix resulting of integrating all its entries. A sufficient condition for this commutation is given by Corollary 2.3 in \cite{TW}. However, it is not necessary and we are not using it in this presentation. Our hypothesis is just (\ref{39}).

\medskip

Then \cite{TW}, its general solution can be written as $ \mathbf x(t)=\Phi(t)\,\mathbf x(0)$, where the initial condition $\mathbf x(0)$ is arbitrary and the fundamental matrix $\Phi(t)$ is given by

\begin{equation}\label{40}
\Phi(t) = \exp\left\{ \int_0^t A(v)\,dv   \right\}\,.
\end{equation}

\medskip

Since we are mainly interested in equations of the type (\ref{20}),  we shall restrict ourselves to the case $n=2$. First of all, let us use the following notation:

\begin{equation}\label{41}
B(t):= \int_0^t A(v)\,dv\,,
\end{equation}
so that (\ref{39}) takes the form:

\begin{equation}\label{42}
B(t)\,B'(t)=B'(t)\,B(t)\,.
\end{equation}

We construct the matrix $B'(t)$ by taken the derivative with respect to $t$ of all entries in $B(t)$. A straightforward integration of (\ref{42}) shows that there exists two non-zero constants $\alpha$ and $\beta$ such that, if we denote by $a_{ij}(t)$ the entries of $A(t)$,

\begin{equation}\label{43}
a_{21}(t)=\alpha\,a_{12}(t)\,, \qquad a_{22}(t)=a_{11}(t)+\beta\,a_{12}(t)\,,
\end{equation}
so that

\begin{equation}\label{44}
B(t)= \left( \begin{array}{cc}  f(t)  &  g(t) \\[2ex]  \alpha\,g(t)  &  f(t)+\beta\, g(t)  \end{array}   \right)\,,
\end{equation}
with

\begin{equation}\label{45}
f(t)=\int_0^ t a_{11}(v)\,dv\,,\qquad g(t)=\int_0^t a_{12}(v)\,dv\,.
\end{equation}

The converse is also true, in the sense that (\ref{43}) implies (\ref{42}).

Then, we may obtain the fundamental matrix (\ref{40}) in the following form:

\begin{equation}\label{46}
\Phi(t)=Q(t)\,\exp\left\{ f(t)+\frac12\, \beta\,g(t)   \right\}\,,
\end{equation}
with

\begin{equation}\label{47}
Q(t)=\left(  \begin{array}{cc}  \cosh\left( \frac12 \,\gamma\,g(t)\right)-\frac \beta \gamma\,
\sinh \left( \frac12 \,\gamma\,g(t) \right)            &  \frac2\gamma \, \sinh \left( \frac12 \,\gamma\,g(t) \right)  \\[2ex]  \frac{2\alpha}\gamma\, \sinh \left( \frac12 \,\gamma\,g(t) \right)  &   \cosh\left( \frac12 \,\gamma\,g(t)\right)+\frac \beta \gamma\,
\sinh \left( \frac12 \,\gamma\,g(t) \right)
 \end{array}     \right)\,.
\end{equation}

Here, $\gamma:= \sqrt{4\alpha+\beta^2}$. Observe that the dependence on $t$ of $Q(t)$ comes solely with $g(t)$  and, hence, of $a_{12}(t)$ only. An explicit expression for $\Phi(t)$ is only possible if we know the primitives for $a_{11}(t)$ and  $a_{12}(t)$. Otherwise, we have to resort to numerical estimations of $f(t)$ and $g(t)$.

\medskip

Let us prove that the fundamental matrix is given by (\ref{46}-\ref{47}). We perform this proof into two steps. First of all, take $\beta=0$ and then, remove this condition.

\smallskip

i.) The case $\beta=0$. Now, $B(t)$ may be written as

\begin{equation}\label{48}
B(t)=f(t)\,I+g(t)\left(\begin{array}{cc}  0 & 1 \\[2ex]
\alpha\,g(t)  &  0      \end{array}   \right)\,,
\end{equation}
where $I$ is the identity matrix. Since any matrix commutes with the identity, by exponentiation we have:

\begin{equation}\label{49}
\Phi(t)=\exp\{f(t)\,I\} \cdot \exp \left[ g(t) \left( \begin{array}{cc} 0 & 1 \\[2ex]  \alpha\,g(t) & 0    \end{array}\right)\right]
\end{equation}

In order to calculate the second exponential in (\ref{49}), we proceed by direct exponentiation of the involved matrix. Then, (\ref{49}) becomes

\begin{equation}\label{50}
\Phi(t)= \exp\{f(t)\} \left(\begin{array}{cc}  1+\frac{g^{2}\alpha}{2}+\frac{g^{4}\alpha^{2}}{24}+\dots & g+\frac{g^{3}\alpha
}{6}+\frac{g^{5}\alpha^{2}}{120}+\dots\\ [2ex]   \alpha\left(g+\frac{g^{3}\alpha}{6}+\frac{g^{5}\alpha^{2}}{120}+\dots\right) &
1+\frac{g^{2}\alpha}{2}+\frac{g^{4}\alpha^{2}}{24}+\dots  \end{array}\right)\,.
\end{equation}

It is obvious that the entries of the matrix in (\ref{50}) are Taylor series corresponding to sinh and cosh centred at the origin. Thus we obtain $Q(t)$ as in  (\ref{47}) for $\beta=0$. 

\smallskip

ii.) The case $\beta\ne 0$. The procedure is essentially the same. In this case, we decompose $B(t)$ as

\begin{equation}\label{51}
B(t)= \left\{ f(t)+\frac\beta 2\,g(t)  \right\}\,I+g(t)\left(\begin{array}{cc} -\beta/2 & 1\\[2ex]  \alpha\,g(t)  &  \beta/2\, g(t)   \end{array}\right)\,,
\end{equation}
so that

\begin{equation}\label{52}
\Phi(t)= \exp\left\{ f(t)+\frac\beta 2\,g(t)  \right\}\cdot \exp\left\{ g(t)\left(\begin{array}{cc} -\beta/2 & 1\\[2ex]  \alpha\,g(t)  &  \beta/2\, g(t)   \end{array}\right)   \right\}\,.
\end{equation}

We exponentiate the matrix in (\ref{52}). After some tedious but straightforward manipulations, we finally obtain $Q(t)$ as in (\ref{51}).

\medskip

Comment: since $A(t)$ is periodic with period $T$, we may obtain the Floquet exponents after (\ref{11}) and (\ref{12}). It is rather straightforward to obtain the following expression, whose real part gives the critical exponents:

\begin{equation}\label{53}
\lambda_\pm=\frac 1T\,\left( f(T)+\frac 12 (\beta\pm\gamma)\,g(T)  \right)\,.
\end{equation}

It is important to remark that this procedure makes sense if (\ref{42}) holds. Otherwise, the eigenvalues of $A(t)$ may differ from the critical exponents, as in the case of the Marcus-Yamabe equation discussed in the previous section.  In addition, in our case the critical exponents also coincide with the eigenvalues of the average of $A(t)$ over a period, which is the matrix given by

\begin{equation}\label{54}
\overline {A(t)}:=\frac 1T \int_0^T A(t)\,dt\,.
\end{equation}

This is an important result that we want to underline here. Also, note that $\gamma$ is complex in general, so that if $g(t)\equiv 0$, then the Floquet coefficients (\ref{53}) are real. This latter result coincides with the obtained in \cite{TW}. 

\medskip

In the sequel, we list some straightforward properties.

\begin{itemize}

\item{If $\gamma^2=4\alpha+\beta^2>0$, see after (\ref{47}), then $Q(t)$ is not bounded. This means that there are  not periodic solutions, even for $A(t)$ periodic. }

\item{If $\gamma^2=4\alpha+\beta^2<0$, then $Q(t)$ is bounded. However, if $f(t)+1/2\,\beta\,g(t)$, or equivalently, $a_{11}(t)+1/2\,\beta \,a_{12}(t)$ has no zeros, one uses (\ref{50}) to see that there are no periodic solutions. Note that the existence of an imaginary part in $\lambda_\pm$ does not imply the periodicity of $\Phi(t)$.}

\item{If $\gamma^2=4\alpha+\beta^2<0$ and $g(t)$ is periodic, we conclude after (\ref{53}) that $Q(t)$ is also periodic. Then, an obvious sufficient condition for the periodicity of solutions is that $f(t)+1/2\,\beta\,g(t)=0$, or equivalently, $a_{11}(t)+1/2\,\beta \,a_{12}(t)=0$ for all values of $t$. This latter conditions implies, in particular, that the critical exponents in (\ref{57}) are purely imaginary, since $\gamma^2<1$.}

\item{If $\gamma^2=4\alpha+\beta^2<0$ and $g(t)$ is periodic with period $T$, then $Q(t)$ is also periodic.}

\item{If $\gamma^2=4\alpha+\beta^2<0$ and, in addition,

\begin{equation}\label{55}
\lim_{t\to\infty}\left(  f(t)+\frac12\,g(t)  \right)=K<\infty\,,
\end{equation}
then, the zero solution (also called the trivial solution), $\mathbf x(t)\equiv \mathbf 0$ all $t$, of (\ref{1}) is stable. }

\item{If $\gamma^2=4\alpha+\beta^2<0$ and, in addition,

\begin{equation}\label{56}
\lim_{t\to\infty}\left(  f(t)+\frac12\,g(t)  \right)=-\infty\,,
\end{equation}
then, the zero solution of (\ref{1}) is a global attractor.}

\item{If $\gamma=0$, then $Q(t)\equiv I$, the identity matrix, for all $t$. There are no periodic or even bounded solutions. }

\end{itemize}

\subsection{An example}

The results of the previous subsection allow us to test the method introduced in Section 3. Now, $A(t)$ is a $2\times 2$ periodic matrix so that relations (\ref{43}) are valid. For the independent entries and parameters in (\ref{43}), we choose:

\begin{equation}\label{57}
a_{11}(t)=-1\,,\qquad a_{12}(t)=2+\sin t\,,\qquad \alpha=-1\,,\qquad \beta=0\,.
\end{equation}

Relations (\ref{57}) fully determine $A(t)$, which is periodic with period $2\pi$. Using (\ref{46}) and (\ref{47}), we have

\begin{equation}\label{62}
\Phi(t)=e^{-t} \left( \begin{array}{cc} \cos(1+2t-\cos t) & \sin(1+2t-\cos t) \\[2ex]  -\sin(1+2t-\cos t) & \cos(1+2t-\cos t) \end{array}\right)\,.
\end{equation}

Using these data, let us write equation (\ref{1}) as a second order linear equation. This gives:

\begin{equation}\label{63}
(2+\sin t)\ddot x(t)+(4-\cos t+2\sin t)\dot x(t)+(10-\cos t+13\sin t+6\sin^2 t+\sin^3 t)x(t)=0\,.
\end{equation}

Next, we determine the critical exponents by our method for $n=3$. The result is a double $\lambda_A=-1-1.9998 i$ (we recall that the exponents are unique modulus $2k\pi i/T$,  $k$ being integer). Using the method given in Section 2, formula (\ref{11}), we obtain $\lambda_{num}=-1\pm 1.09\,10^{-6} i$. Here, it is possible to obtain the exact value, which gives $\lambda_e=-1-2i$. The conclusion is that our method gives far more precision than the standard numerical method.  

Finally, let us integrate (\ref{63}) using our method, approaching coefficients by their nearest rational number and use trigonometric relations. We have the following approximation for the solution:

\begin{eqnarray}\label{64}
x_A(t)=e^{-t} \left(  \frac 18\,i+\frac{8}{17}\,\cos  t-\frac 78 \,i\,\cos (2t)+\frac 12\,\cos (3t)  +\frac 18\,i\,\cos(4x)      -\frac1{42}\,\cos(5t)                          \right.  \nonumber\\[2ex]  \left.-\frac{10}{19}\,i\,\sin t -\frac{26}{30}\,\sin(2t) -\frac12\,i\,\sin(3t)+\frac18\,\sin(4t)+\frac1{42} \,i\,\sin(5t) +\dots
\right)\,.
\end{eqnarray} 

Let us obtain the first terms of the exact solution by using the initial conditions  $x(0)=x_A(0)$ and $\dot x(0)=\dot x_A(0)$ and expanding the entries of $\Phi(t)$ in Fourier series. The result is

\begin{eqnarray}\label{65}
x_e(t)=e^{-t} \left(  \frac 18\,i+\frac{7}{15}\,\cos  t-\frac 78 \,i\,\cos (2t)+\frac 12\,\cos (3t)  +\frac 18\,i\,\cos(4x)      -\frac1{45}\,\cos(5t)                          \right.  \nonumber\\[2ex]  \left.-\frac{6}{7}\,i\,\sin t -\frac{26}{30}\,\sin(2t) -\frac12\,i\,\sin(3t)+\frac18\,\sin(4t)+\frac1{45} \,i\,\sin(5t) +\dots
\right)\,.
\end{eqnarray} 

The coincidence between both results is high showing the remarkable accuracy of our method.

\subsection{A second example: The Marcus-Yamabe equation}

As a second example, let us consider the Marcus-Yamabe system, which is of second order. This has been given as a counter-example of a periodic system such that the eigenvalues of $A(t)$ are constant, i.e., independent of $t$, equal to $(-1\pm \sqrt 7\,i)/4$ and yet being the zero solution not stable   \cite{VER,FAR}. The Marcus Yamabe system is of the form (\ref{1}), with matrix $A(t)$ given by 

\begin{equation}\label{66}
A(t)=\left( \begin{array}{cc}  -1+\frac32\,\cos^2 t  &  1-\frac32 \,\cos t\, \sin t   \\[2ex] -1-\frac32\,\cos t\,\sin t &  -1+\frac32\,\sin^2 t  \end{array}   \right)\,.
\end{equation}

System (\ref{64}) has two linearly independent solutions of the form:

\begin{equation}\label{67}
\mathbf x_1(t)= \left(\begin{array}{c}  -\cos t\\[2ex] \sin t  \end{array}  \right)\, e^{t/2}\,, \qquad \mathbf x_2(t)= \left( \begin{array}{c} \sin t \\[2ex] \cos t   \end{array}   \right)\, e^{-t}\,.
\end{equation}

It is easy to write the associated differential equation of the Marcus-Yamabe system, the Marcus-Yamabe equation. This is

\begin{equation}\label{68}
(8-6\sin 2t)\ddot x(t)+(4+12\cos 2t-3 \sin 2t) \dot x(t)+(-5+3\cos 2t+9\sin 2t)x(t)=0\,.
\end{equation}

This equation does not have singular points. Its Floquet characteristic coefficients are $-1$ and $1/2$, see (\ref{67}). If we apply the method introduced in Section 3 with $n=3$, we obtain the same critical exponents and the following basis of the space of solutions:

\begin{equation}\label{69}
x_1(t)=e^{-t}\,\sin t\,,\qquad x_2(t)=e^{t/2}\,\cos t\,,
\end{equation} 
with coincide with the exact solution as we see from (\ref{67}). The conclusion is that our approximate method has yield to the exact solution with $n=3$.

\section{Concluding remarks}

The Hill-Harmonic Balance method has been designed in order to obtain the Floquet characteristic exponents for linear differential equations and systems with periodic coefficients. These exponents are solutions of an algebraic equation of degree $2n+1$, where $n$ is the order of   a Fourier polynomial that it is used in order to obtain an approximate analytic solution for the equation.  Since in general, $2n+1$ is much larger than the order of the equation, that in many practical cases is two, we need an efficient method so as to choose the Floquet exponents among the solutions of the algebraic equation. They are determined by a variational method.  It is precisely  the use of this variational principle that determines the Floquet exponents as its critical points, which makes our procedure different from other discussed in the literature. These tools are easy to implement for practical applications. We obtain an excellent precision in function of the number of harmonics used. 

We have compared our results with the exact results known for the Mathieu equation. They show a good accuracy even if we just take the first two nodes (up to $n=2$) in the Fourier series. The precision obtained for $n=3$ is excellent. We also have compared our results with those obtained with the standard method described in Section 2. The conclusion is that we obtain better results with little effort and negligible computational time. 

We have also analysed the particular situation where a $2\times 2$ periodic matrix $A(t)$ commutes with its integral. An example shows the high precision of the method also in this situation. The critical exponents are the eigenvalues of the average of the matrix $A(t)$ over a period. 

We have also tested our results in another differential equation known as the Marcus-Yamabe equation. In this case a three nodes approximation, $n=3$, gives the exact result.

In conclusion: this is a method to obtain Floquet characteristic coefficients which is simple, efficient and with an excellent precision as shown in the testing examples.   Although we have not proposed an explicit formula to evaluate the error, once we have determined critical exponents and approximate solutions, formula (\ref{24}) may serve to test the accuracy of a given solution.

\section*{Acknowledgements}

We acknowledge partial financial support to the Spanish Ministry of Economy grant  MTM2014-57129-C2-1-P, to the Junta de Castilla y Le\'on grant VA057U16 as well as to the National University of Rosario (Argentina) grant ING19/i402.

\end{document}